%% file: FiniteSupport_n_alpha.tex
\documentclass[4pt, oneside]{article}

\usepackage[OT2, T1]{fontenc}
\usepackage[russian, english]{babel}
\usepackage{setspace}

\usepackage{enumerate}
\input{Preamble}

\newcommand{\TPr}{\textbf{\foreignlanguage{russian}{TS}C}}
\newcommand{\al}{\alpha}
\newcommand{\be}{\beta}
\newcommand{\FLE}{\mathbb{F}}
\newcommand{\MNF}[2]{\la n_{#1}^{\al_{#1}} \ra \top \land \ldots \land \la n_{#2}^{\al_{#2}} \ra \top }

\newcommand{\Lim}{\text{Lim}}

\def\le{{\ell}}

\def\I{{\ensuremath{\mathcal{J}}}\xspace}

\begin{document}

\title{A Finitely Supported Frame for the Turing Schmerl Calculus}
\author{Eduardo Hermo Reyes \footnote{\href{mailto:ehermo.reyes@ub.edu}{ehermo.reyes@ub.edu} } \\ {\small University of Barcelona}\\ {\small Department of Philosophy}  }
\maketitle

\begin{abstract}
In \cite{HermoJoosten:2017:TSCArith} we introduced the propositional modal logic $\TPr$ (which stands for Turing Schmerl Calculus) which adequately describes the provable interrelations between different kinds of Turing progressions. In \cite{HermoJoosten:2018:TSCSem_omegaLambda} we defined a model $\mathcal{J}$ which is proven to be a universal model for $\TPr$ based on the intensively studied Ignatiev's universal model for the closed fragment of  $\textbf{GLP}$ (G\"odel L\"ob's polymodal provability logic). In the current paper we present a new universal frame $\mathcal{H}$, which is a slight modification of $\mathcal{J}$, and whose domain allows for a modal definability of each of its worlds.
\end{abstract}

\section{Introduction}
Turing progressions arise by iteratedly adding consistency statements to a base theory. Different notions of consistency give rise to different Turing progressions. In \cite{HermoJoosten:2017:TSCArith}, we introduced the system $\TPr$ that generates exactly all relations that hold between these different Turing progressions given a particular set of natural consistency notions. The system was proven to be arithmetically sound and complete for a natural interpretation, named the \emph{Formalized Turing progressions} (FTP) interpretation. In \cite{HermoJoosten:2018:TSCSem_omegaLambda}, we defined a model $\mathcal{J}$ which is proven to be a universal model for $\TPr$. The model $\mathcal{J}$  is a slight modification of the intensively studied $\mathcal{I}$ : Ignatiev's universal model for the closed fragment of G\"odel L\"ob's polymodal provability logic $\textbf{GLP}$.\\

In the current paper we introduce a new universal model for $\TPr$ based on a minor modification of $\mathcal{J}$. Ignatiev's model and $\mathcal{J}$ share the same domain built from specials sequences of ordinals. For this new frame, we shall consider only sequences with finite support. We will prove the completeness of $\TPr$ w.r.t. this model and as a consequence, we shall see that every world in the frame is modally definable.

\section{The logic $\TPr$}

In this section we introduce the logic $\TPr$ whose main goal is to express valid relations that hold between the corresponding Turing progressions. $\TPr$ is built-up from a positive propositional modal signature  using \emph{ordinal modalities}. Let $\Lambda$ be a fixed recursive ordinal throughout the paper with some properties as specified in Remark \ref{remark:WhyEpsilonNumbers}. By ordinal modalities we denote modalities of the form $\la n^\al  \ra$ where $\al \in \Lambda$ for some fixed ordinal $\Lambda$ and $n \in \omega$ (named \emph{exponent} and \emph{base}, respectively). The set of formulas in this language is defined as follows:

\begin{definition}By $\FLE$ we denote the smallest set such that:
\begin{enumerate}[i)]
\itemsep0em
\item $\top \in \FLE$;
\item If $\varphi, \, \psi \in \FLE\Rightarrow (\varphi \wedge \psi) \in  \FLE$;
\item if $\varphi \in \FLE, \ n < \omega$ and $\al < \Lambda \Rightarrow \la n^\al \ra \varphi \in \FLE$.
\end{enumerate}
\end{definition}

For any formula $\psi$ in this signature, we define the set of base elements occurring in $\psi$. That is:

\begin{definition}
The set of base elements occurring in any modality of a formula $\psi \in \FLE$ is denoted by ${\sf N\text{-}mod} (\psi)$. We recursively define ${\sf N\text{-}mod}$ as follows:
\begin{enumerate}[i)]
\item ${\sf N\text{-}mod}(\top) = \emptyset$;
\item ${\sf N\text{-}mod}(\varphi \wedge \psi) = {\sf N\text{-}mod}(\varphi) \cup {\sf N\text{-}mod}(\psi)$;
\item ${\sf N\text{-}mod}(\la n^\al \ra \psi) = \{ n \} \cup {\sf N\text{-}mod}(\psi)$.	
\end{enumerate}
\end{definition}

In order to give the axioms of our system, we need to consider a kind of special formulas named \emph{monomial normal forms}. Monomial normal forms are conjunctions of monomials with an additional condition on the occurring exponents. In order to formulate this condition we first need to define the \emph{hyper-exponential} as studied in \cite{FernandezJoosten:2012:Hyperations}.

\begin{definition}
For every $n \in \omega$ the \emph{hyper-exponential} functions $e^n: \text{On} \rightarrow \text{On}$ are recursively defined as follows: $e^0 $ is the identity function, $e^1: \al \mapsto -1+\omega^\al$ and $e^{n+m}= e^n \circ e^m$.
\end{definition}

We will use $e$ to denote $e^1$. Note that for $\alpha$ not equal to zero we have that $e(\alpha)$ coincides with the regular ordinal exponentiation with base $\omega$; that is, $\alpha \mapsto \omega^\alpha$. However, it turns out that hyper-exponentials have the nicer algebraic properties in the context of provability logics.

\begin{definition}
The set of formulas in \emph{monomial normal form}, {\sf MNF}, is inductively defined as follows:
\begin{enumerate}[i)]
\item $\top \in {\sf MNF}$;
\item $\la n^\al \ra \top \in {\sf M}$, for any $n \, {<} \, \omega$ and $\al \, {<} \, \Lambda$;
\item \begin{tabbing}
	  if \hspace{0.2cm} \= $a) \ $ \= $\la n_0^{\al_0} \ra \top \wedge \ldots \wedge \la n_k^{\al_k} \ra \top \in {\sf MNF}$; \\
 	  \	\\	     
	     \> $b) \ $ \> \ $n < n_0$; \\
	  \ \\
	     \> $c) \ $ \> \ $\al$ of the form $e^{n_0 - n}(\al_0) \cdot (2 + \delta)$ for some $\delta < \Lambda$,
      \end{tabbing}
      \vspace{0.25cm}
then $\la n^\al \ra \top \wedge \la n_0^{\al_0} \ra \top \wedge \ldots \wedge \la n_k^{\al_k} \rangle \top \in {\sf MNF}.$
\end{enumerate}
\end{definition}

The derivable objects of $\TPr$ are \emph{sequents} i.e. expressions of the form $\varphi \vdash \psi$ where $\varphi, \, \psi \in \mathbb{F}$. We will use the following notation: by $\varphi \equiv \psi$ we will denote that both $\varphi \vdash \psi$ and $\psi \vdash \varphi$ are derivable. Also, by convention we take that for any $n$, $\la n^0 \ra \varphi$ is just $\varphi$.

\begin{definition} 
$\TPr$ is given by the following set of axioms and rules:\\

Axioms:
\begin{enumerate}

\item $\varphi \vdash \varphi, \ \ \ \varphi \vdash \top$;	

\item $\varphi \wedge \psi \vdash \varphi,  \ \ \ \varphi \wedge \psi \vdash \psi$;	\label{conjel}

\item Monotonicity axioms: $\la n^\al \ra \varphi \vdash \la n^\be \ra \varphi$, \ \ \ for $\be \, {<} \, \al$;	\label{mon}
		
\item Co-additivity axioms: $\la n^{\be + \al} \ra \varphi \equiv \la n^\al \ra \la n^\be \ra \varphi$;	 \label{coadditive}

\item Reduction axioms: $\la (n + m)^\al \ra \varphi \vdash \la n^{e^m (\al)} \ra \varphi$;	\label{omega}

\item Schmerl axioms: 
\[
\la n^\al \ra \big( \, \la n_0^{\al_0} \ra \top \ \wedge \ \psi \, \big) \equiv \la n^{e^{n_0 - n} (\al_0) \cdot (1 + \al)}  \ra \top \land \la n_0^{\al_0} \ra \top \ \land \ \psi
\] 
for $n\, {<} \, n_0$ and $\la n_0^{\al_0} \ra \top \ \wedge \ \psi \in {\sf MNF}$. \label{MS1}

\end{enumerate}
\

Rules:
\begin{enumerate}
\item If $\varphi \vdash {\psi}$ and $\varphi\vdash{\chi}$, then $\varphi\vdash {\psi}\wedge{\chi}$;	\label{r:1}
\item If $\varphi\vdash {\psi}$ and ${\psi}\vdash{\chi}$, then $\varphi\vdash{\chi}$;	\label{r:2}
\item If $\varphi\vdash {\psi}$, then $\la n^\al \ra \varphi\vdash \la n^\al \ra {\psi}$ ;	\label{r:3}
\item If $\varphi \vdash \psi$, then $\la n^\al \ra \varphi \, \land \, \la m^{\be + 1} \ra \psi \, \vdash \, \la n^\al \ra \big( \, \varphi \, \land \, \la m^{\be+ 1} \ra \psi \, \big)$ \ for $n \, {>} \, m$. \label{r:4}
\end{enumerate}
\end{definition}

It is worth mentioning the special character of Axioms (\ref{omega}) and (\ref{MS1}) since both axioms are modal formulations of principles related to Schmerl's fine structure theorem, also known as \emph{Schmerl's formulas} (see \cite{Schmerl:1978:FineStructure} and \cite{Beklemishev:2003:AnalysisIteratedReflection}).

\begin{remark}\label{remark:WhyEpsilonNumbers}
As we see in the axioms of our logic, they only make sense if the ordinals occuring in them are available. Recall that $\Lambda$ is fixed to be a recursive ordinal all through the paper. Moreover, some usable closure conditions on $\Lambda$ naturally suggest themselves. Since it suffices to require that for $n<\omega$ that $\alpha,\beta < \Lambda \ \Rightarrow \ \alpha + e^n(\beta) < \Lambda$, we shall for the remainder assume that $\Lambda$ is an $\varepsilon$-number, that is, a positive fixpoint of $e$ whence $e(\Lambda)=\Lambda=\omega^\Lambda$.
\end{remark}

In \cite{HermoJoosten:2017:TSCArith}, the authors proved that for any formula $\varphi$, there is a unique equivalent $\psi$ in monomial normal form. 

\begin{theorem} \label{MNFT}
For every formula $\varphi$ there is a unique $\psi \, {\in} \, {\sf MNF}$ such that $\varphi \equiv \psi$.
\end{theorem}

In virtue of the Reduction axioms, a formula $\psi \in {\sf MNF}$ may bear implicit information on monomials $\la n^\al \ra \top$ for $n \not \in {\sf N\text{-}mod}(\psi)$. The next definition is made to retrieve this information.

\begin{definition} \label{projection}
Let $\psi := \MNF{0}{k} \in {\sf MNF}$. By $\pi_{n_i} (\psi)$ we denote the corresponding exponent $\al_i$. Moreover, for $m \not \in {\sf N\text{-}mod}(\psi)$, with $n_k > m$, $\pi_m ( \psi)$ is set to be $e\big(\, \pi_{m+1} (\psi) \, \big)$ and for $m' > n_k$, $\pi_{m'} (\psi)$ is defined to be $0$.
\end{definition}

Another important result in \cite{HermoJoosten:2017:TSCArith} is the decidability of the sequents entailing formulas in monomial normal form:

\begin{theorem} 
\label{CharacterizatioMNFDeriv}
For any $\psi_0, \ \psi_1 \in {\sf MNF}$, where $\psi_0 := \MNF{0}{k}$ and $\psi_1 := \la m_0^{\be_0} \ra \top \ \land \ \ldots \ \land \ \la m_j^{\be_j} \ra \top$. We have that $\psi_0 \vdash \psi_1$ iff for any $n < \omega$, $\pi_n (\psi_0) \geq \pi_n (\psi_1)$.
\end{theorem}

\section{A universal model for $\TPr$}

In \cite{HermoJoosten:2018:TSCSem_omegaLambda} we defined a modal model $\mathcal{J}$ which is universal for $\TPr$. That is, any derivable sequent holds everywhere in the model whereas any non-derivable sequent is refuted somewhere in the model. This section summarizes the results in that paper.\\

As Ignatiev's frame, this model is based on specials sequences of ordinals. In order to define them, we need the following central definition.

\begin{definition}
We define \emph{ordinal logarithm} as $\le (0):= 0$ and $\le (\alpha + \omega^\beta) := \beta$.
\end{definition}

With this last definition we are now ready to introduce the set of worlds of our frame.

\begin{definition}
By $\text{Ig}^\omega $ we denote the set of \emph{$\ell$-sequences} or \emph{Ignatiev sequences}. That is, the set of sequences $x :=  \la x_0, x_1, x_2, \ldots \ra$ where for $i \, {<} \, \omega$, $x_{i+1}\leq \le(x_i)$. 
\end{definition}

Next, we can define our frame, which is a minor variation of Ignatiev's frame.

\begin{definition} 
$\mathcal{J}_\Lambda := \la I, \{ R_n \}_{n < \omega}\ra$ is defined as follows:
$$I := \{ x \, {\in} \, \text{Ig}^\omega : x_i \, {<} \, \Lambda \ \text{for } i \, {<} \, \omega \}$$ and
$$x R_n y :\Leftrightarrow (\forall \, m \leq n \ x_m \, {>} \, y_m \ \wedge \ \forall \, i \, {>} \, n \ x_i \geq y_i).$$
\end{definition}

Since $\Lambda$ is a fixed ordinal along the paper, from now on we suppress the subindex $\Lambda$. \\

We define the auxiliary relations $R_n^\alpha$ for any $n < \omega$ and $\alpha < \Lambda$. The idea is that the $R_n^\alpha$ will model the $\la n^\alpha \ra$ modality.

\begin{definition} \label{JOrdRelation}
Given $x, y \in I$ and $R_n$ on $I$, we recursively define $x R_n^\alpha y$ as follows:
\begin{enumerate}
\item 
$x R_n^0 y \ \ :\Leftrightarrow \ \ x=y$;

\item 
$x R_n^{1+\alpha} y \ :\Leftrightarrow \ \forall \, \beta  {<} 1{+}\alpha  \ \exists z \ \big( xR_nz \ \wedge \ z R_n^\beta y\big)$.
\end{enumerate}
\end{definition}

With the the auxiliary relations $R_n^\alpha$, we can define the validity of a formula $\varphi$ in a point $x \in I$. 

\begin{definition} \label{FormulaTrueInPoint}
Let $x \in I$ and $\varphi \in \FLE$. By $x \Vdash \varphi$ we denote the validity of $\varphi$ in $x$ that is recursively defined as follows:
\begin{itemize}
\item $x \Vdash \top$ for all $x \in I$; 
\item $x \Vdash \varphi \wedge \psi$ iff $x \Vdash \varphi$ and $x \Vdash \psi$;
\item $x \Vdash \la n^\alpha \ra \varphi$ iff there is $y \in I$, $x R_n^\alpha y$ and $y \Vdash \varphi$.
\end{itemize}
\end{definition}

The following theorem establishes the completeness of $\TPr$.

\begin{theorem} \label{JCompleteness}
For any $\varphi,\, \psi \in \FLE$, we have that:
\[
\varphi \vdash \psi \ \Longleftrightarrow \ \forall \, x \in I \ \Big( \, \mathcal{J}, x \Vdash \varphi \, \Longrightarrow \, \mathcal{J}, x \Vdash \psi \, \Big).
\]
\end{theorem}

\subsection{Ignatiev sequences and MNF's}

In this subsection we shall define a translation between formulas in formulas in monomial normal form and Ignatiev sequences. 

\begin{definition}\label{TranslationMNFI}
Let $\psi:= \MNF{0}{k} \in {\sf MNF}$. By $x_\psi$ we denote the sequence $\la \pi_i (\psi) \ra_{i < \omega}$.
\end{definition}

In virtue of Definition \ref{projection}, we can observe that for $\psi \in {\sf MNF}$, we have that $x_\psi \in I$.\\

The following results establish the relation between this translation and the validity in $\mathcal{J}$:

\begin{lemma} \label{trueInTrans}
For any $\varphi \in {\sf MNF}, \ \mathcal{J}, x_\varphi \Vdash \varphi$.
\end{lemma}

\begin{lemma} \label{failsInTrans}
Given $\varphi, \, \psi \in {\sf MNF}$, if $\varphi \not \vdash \psi$ then $\mathcal{J}, x_\varphi \Vdash \varphi$ and $\mathcal{J},x_\varphi \not \Vdash \psi$.
\end{lemma}

\section{A finitely supported variation on $\mathcal{J}$}

For $\Lambda > \varepsilon_0$ it might happen that a $\ell$-sequence $x$ never reach zero, getting stabilized at some $\varepsilon$-number. For instance, the sequence $\la \varepsilon_0,\, \varepsilon_0,\, \varepsilon_0,\, \ldots \, \ra \in \text{Ig}^\omega $ but since $\varepsilon_0$ is a fixpoint of ordinal exponentiation with base $\omega$, that is, $\omega^{\varepsilon_0} = \varepsilon_0$, this sequence will not converge to zero.\\

Given a $\ell$-sequence $x$, if all but finitely many of its elements are zero, we will write $\la x_0 , \ldots, x_n, \vec 0\ra$ to denote such $\le$-sequence or even simply $\la x_0 , \ldots, x_n\ra$ whenever $x_{n+1} =0$. We shall refer to such sequences as \emph{finitely supported sequences} or \emph{sequences with finite support}. In \cite{HermoJoosten:2018:TSCSem_omegaLambda}, we showed that only sequences with finite support are modally definable. 

The purpose of this new universal frame is to give a model where every world is modally definable. Thus, the set of worlds of our frame is going to be built from these finitely supported sequences whose ordinals are less than $\Lambda$. Thereby, the domain of our frame is a subset of $I$. With respect to the relations, we maintain the same definition that we introduced for the relations in $\mathcal{J}$. 

\begin{definition}
$\mathcal{H}_\Lambda := \la H, \{ S_n \}_{n < \omega}\ra$ is defined as follows:
$$H := \{ x \, {\in} \, \text{Ig}^\omega : \forall \, i \, {<} \, \omega \ x_i \, {<} \, \Lambda \ \land \ \exists \, j \, {<} \, \omega \ x_j = 0 \ \}$$ and
$$x S_n y :\Leftrightarrow (\forall \, m \leq n \ x_m \, {>} \, y_m \ \wedge \ \forall \, i \, {>} \, n \ x_i \geq y_i).$$
\end{definition}

As before, since $\Lambda$ is a fixed ordinal along the paper, from now on we suppress the subindex $\Lambda$.

\begin{remark} \label{rmk1}
Since we have that $H \subset I$ and we preserve the same definition for the $S_n$ relations, clearly for any $x, \, y \in H$, we have that $x S_n y \ \Longleftrightarrow \ x R_n y$. Moreover, if $x \in H$ and $x R_n y$ for some $y \in I$, then $y \in H$ and $x S_n y$.
\end{remark}

We can make some simple observations on the relations $S_n$.

\begin{lemma}\label{theorem:G_basicPropertiesRnRelations}\ 
\begin{enumerate}
\item
Each $S_n$ for $n\in \omega$ is
\emph{transitive}: $xS_ny \ yS_nz \ \Rightarrow \ xSz$;

\item\label{theorem:G_basicPropertiesRnRelations:Item:Noetherian}
Each $S_n$ for $n\in \omega$ is
\emph{Noetherian}: each non-empty $X\subseteq H$ has an $S_n$-maximal element $y\in X$, i.e., $\forall\, x {\in} X\ \neg yS_nx$;

\item
The relations $S_n$ are \emph{monotone} in $n$ in the sense that: $xS_ny \Rightarrow xS_my$ whenever $n>m$.

\end{enumerate}
\end{lemma}

As we did for $\mathcal{J}$, we need to define the auxiliary relations $S_n^\alpha$ for any $n < \omega$ and $\alpha < \Lambda$ whose purpose is to model the $\la n^\alpha \ra$ modality.

\begin{definition} \label{OrdRelation}
Given $x, y \in H$ and $S_n$ on $H$, we recursively define $x S_n^\alpha y$ as follows:
\begin{enumerate}
\item 
$x S_n^0 y \ \ :\Leftrightarrow \ \ x=y$;

\item 
$x S_n^{1+\alpha} y \ :\Leftrightarrow \ \forall \, \beta  {<} 1{+}\alpha  \ \exists z \ \big( xS_nz \ \wedge \ z S_n^\beta y\big)$.
\end{enumerate}
\end{definition}

With the the auxiliary relations $S_n^\alpha$, we give the following definition for a formula $\varphi$ being true in a point $x$ of $\mathcal{H}$. 

\begin{definition} \label{FormulaTrueInPoint}
Let $x \in H$ and $\varphi \in \FLE$. By $x \Vdash \varphi$ we denote the validity of $\varphi$ in $x$ that is recursively defined as follows:
\begin{itemize}
\item $x \Vdash \top$ for all $x \in H$; 
\item $x \Vdash \varphi \wedge \psi$ iff $x \Vdash \varphi$ and $x \Vdash \psi$;
\item $x \Vdash \la n^\alpha \ra \varphi$ iff there is $y \in H$, $x S_n^\alpha y$ and $y \Vdash \varphi$.
\end{itemize}
\end{definition}

The following result was proven in \cite{HermoJoosten:2018:TSCSem_omegaLambda} for the relations $R_n^{\al+1}$ defined over $\mathcal{J}$. The same result holds if we restrict ourselves to elements in $H$ and thus, it can be easily applied to the $S_n^{\al +1}$ relations.

\begin{proposition} Given $x,\, y \in H$, $n < \omega$ and $\al < \Lambda$:
\begin{enumerate}
\item $x R_n^{\al + 1} y \ \Longleftrightarrow \ \exists z \in H \ \big(xR_nz \ \wedge \ z R_n^\al y \big)$;\label{SucRelation1}
\item $x S_n^{\al + 1} y \ \Longleftrightarrow \ \exists z \ \big(xS_nz \ \wedge \ z S_n^\al y \big)$. \label{SucRelation2}
\end{enumerate}
\label{SucRelation}
\end{proposition}
\begin{proof}
We proceed by making a case distinction on $\al$. For $\al=0$ both items are straightforward. For $\al > 0$, we start by proving left-to-right implication of Item \ref{SucRelation1}. Assume $x R_n^{\al + 1} y$. By Definition \ref{JOrdRelation}, $\forall \be < \al + 1 \, \exists z \ \big( x R_n z \ \wedge \ z R_n^\be y \big)$. Since $x \in H$ and $x R_n z$, then $z \in H$. Therefore, in particular we have that $\exists z \in H \ \big( x R_n z \ \wedge \ z R_n^\al y \big)$. For the other direction, assume $\exists z \in H \ \big( x R_n z \ \wedge \ z R_n^\al y \big)$. Thus, $\exists z \in H \ \big( x R_n z \ \wedge \ z R_n^\al y \big)$ and $\forall \be < \al \, \exists u \ \big( z R_n u \ \wedge \ u R_n^\be y \big)$. Hence, we get that $\forall \be < \al + 1 \, \exists z \ \big( x R_n z \ \wedge \ z R_n^\be y \big)$ i.e. $x R_n^{\al + 1} y$. The proof of Item \ref{SucRelation2} is analogous.
\end{proof}

The following proposition establishes the preservation between the $S_n^\al$ relations of $\mathcal{H}$ and the $R_n^\al$ relations of $\mathcal{J}$ when restricted to the elements in $H$.

\begin{proposition} \label{preserveRelation}
For any $x, \, y \in H$, $n < \omega$ and $\al < \Lambda$:
\[
x S_n^\al y \ \Longleftrightarrow \ x R_n^\al y.
\]
\end{proposition}
\begin{proof}
We proceed by induction on $\al$. For $\al = 0$ is trivial. \textbf{For} $\boldsymbol{\al = \be + 1}$, assume $ x S_n^{\be + 1} y$ and reason as follows:
\begin{tabbing}
$x S_n^{\be + 1} y$ \= $\Longleftrightarrow \ \exists z \ \big(xS_nz \ \wedge \ z S_n^\be y \big)$   by Lemma \ref{SucRelation}, Item \ref{SucRelation1};\\[0.30cm]
\ \> $\Longleftrightarrow \ \exists z \in H \ \big(xR_nz \ \wedge \ z R_n^\be y \big)$ by Remark \ref{rmk1} and the I.H.;\\[0.30cm]
\ \> $\Longleftrightarrow \ x S_n^{\be + 1} y$  by Lemma \ref{SucRelation}, Item \ref{SucRelation1}.
\end{tabbing}
\textbf{For} $\boldsymbol{\al = \lambda \in \Lim}$ we reason analogously applying the I.H. and Remark \ref{rmk1}.
\end{proof}

\section{Preservation of validity}

To prove the completeness of $\TPr$ w.r.t. $\mathcal{H}$, first we shall prove the preservation of validity between frames $\mathcal{H}$ and $\mathcal{J}$ when restricted to $\ell$-sequences with finite support. 

\begin{theorem} \label{PreservationVal}
For any $x \in H$ and $\varphi \in \FLE$,
\[
\mathcal{J}, x \Vdash \varphi \ \Longleftrightarrow \ \mathcal{H}, x \Vdash \varphi.
\]
\end{theorem}

\begin{proof}
Proof goes by induction on $\varphi$. Base and conjunctive cases are straightforward. For the case $\varphi := \la n^\al \ra \psi$, assume $\mathcal{J}, x \Vdash \la n^\al \ra \psi$. Thus, there is $y \in I$ such that $x R_n^\al y$ and $\mathcal{J}, y \Vdash \psi$. By Remark \ref{rmk1}, we know $y \in H$ and by Proposition \ref{preserveRelation} we have that $x S_n^\al y$. This combine with the I.H. gives us that there is $y \in H$ such that $x S_n^\al y$ and $\mathcal{H}, y \Vdash \psi$, that is, $\mathcal{H}, x \Vdash \la n^\al \ra \psi$. For the other implication, suppose that $\mathcal{H}, x \Vdash \la n^\al \ra \psi$ i.e. there is $y \in H$ such that $x S_n^\al y$ and $\mathcal{H}, y \Vdash \psi$. By Remark \ref{rmk1}, we know $H \subset I$. Thus, by Proposition \ref{preserveRelation} and the I.H. we obtain that $y \in I$ such that $x R_n^\al y$ and $\mathcal{J}, y \Vdash \psi$ and so, $\mathcal{J}, x \Vdash \la n^\al \ra \psi$.  
\end{proof}

Now we can show the completeness of our system w.r.t. $\mathcal{H}$.

\begin{theorem}
For any $\varphi,\, \psi \in \FLE$, we have that:
\[
\varphi \vdash \psi \ \Longleftrightarrow \ \forall \, x \in H \ \Big( \, \mathcal{H}, x \Vdash \varphi \, \Longrightarrow \, \mathcal{H}, x \Vdash \psi \, \Big).
\]
\end{theorem}
\begin{proof}
For the left-to-right implication, assume $\varphi \vdash \psi$. Thus, by Theorem \ref{JCompleteness}: 
\begin{equation} \label{Eq1}
\forall \, x \in I \ \Big( \, \mathcal{J}, x \Vdash \varphi \, \Longrightarrow \, \mathcal{J}, x \Vdash \psi \, \Big).
\end{equation}

Assume $\mathcal{H}, x \Vdash \varphi$ for some $x \in H$. Then, by Theorem \ref{PreservationVal}, $\mathcal{J}, x \Vdash \varphi $ and by (\ref{Eq1}),
$\mathcal{J}, x \Vdash \psi$. Hence, by Theorem \ref{PreservationVal}, $\mathcal{H}, x \Vdash \psi$. \\

For the other implication, reasoning by contrapositive, assume $\varphi \not \vdash \psi$. Thus, by Lemma \ref{failsInTrans}, we have the following:
\begin{equation} \label{Eq2}
\mathcal{J}, x_\varphi \Vdash \varphi \ \& \ \mathcal{J}, x_\varphi \not \Vdash \psi.
\end{equation}

We can observe that by definitions \ref{TranslationMNFI} and \ref{projection}, $x_\varphi$ is a $\ell$-sequence with finite support. Thus, $x_\varphi \in H$. Hence, by means of Theorem \ref{PreservationVal} together with (\ref{Eq2}), we obtain that $\mathcal{H}, x_\varphi \Vdash \varphi$ but $\mathcal{H}, x_\varphi \not \Vdash \psi$ which concludes the proof.
\end{proof}

\section{Modal definability}

In this section we shall see an application of the completeness of $\TPr$ w.r.t $\mathcal{H}$. As pointed by some results established in \cite{HermoJoosten:2018:TSCSem_omegaLambda}, finitely supported sequences are modally definable, which makes every world in $\mathcal{H}$ modally definable in the sense of Theorem \ref{ModalDefinability}. \\

First, we will check some facts that relate monomial normal forms and $\ell$-sequences with finite support.

\begin{lemma} \label{GtoMNF}
Let $x \in H$. There is a unique $\psi \in {\sf MNF}$ such that $x_\psi = x$.
\end{lemma}
\begin{proof}
Let $x := \la x_0, \ldots, x_m, 0\ra \in H$ and consider the map $M : H \rightarrow {\sf MNF}$ defined as follows:
\begin{enumerate}
\item $M(\la 0 \ra) = \top$;
\item $M(\la \al , 0 \ra) = \la 0^{\al} \ra \top$;
\item 
$M(\la \al_0, \ldots, \al_{j+1}, 0 \ra) = {\small \begin{cases} M(\la \al_1, \ldots, \al_{j+1}, 0 \ra) & \text{if } \al_0 = e(\al_1);\\[0.20cm] \la 0^{\al_0} \ra \top \ \wedge \ M(\la \al_1, \ldots, \al_{j+1}, 0 \ra) & \text{if } \al_0 > e(\al_1). \end{cases}}$
\end{enumerate}
By Definition \ref{TranslationMNFI}, we can easily check that $x_{M(x)} =x$.
\end{proof}

With the following auxiliary definition we are ready to introduce the main theorem of this section about the modal definability of the worlds in $\mathcal{H}$.

\begin{definition}
Let $x \in H$. By $x^\downarrow$ we denote the subset of $H$ such that $y \in x^\downarrow$ iff $y_j < x_j$ for some $j < \omega$.
\end{definition}

\begin{theorem} \label{ModalDefinability}
For any $x \in H$ there is $\varphi \in \FLE$ such that: 
\[
\mathcal{H}, x \Vdash \varphi \ \ \& \ \ \forall \, y \in x^\downarrow, \ \mathcal{H}, y \not \Vdash \varphi.
\]
\end{theorem}

\begin{proof}
By a combination of Lemma \ref{GtoMNF} and Lemma \ref{trueInTrans} we obtain that for every $x \in H$, we have that $\mathcal{H}, x \Vdash M(x)$. Moreover, by Theorem \ref{CharacterizatioMNFDeriv}, we have that for any $y \in x^\downarrow$, $M(y) \not \vdash M(x)$. Hence, by Lemma \ref{failsInTrans} and  Theorem \ref{PreservationVal}, $\mathcal{H}, y_{M(y)} \not \Vdash M(x)$ i.e. $\mathcal{H}, y \not \Vdash M(x)$.
 
\end{proof}

\nocite{*}
\bibliographystyle{plain}
\bibliography{ref1.bib}

\end{document}

%% file: Preamble.tex
\usepackage{amsmath,amssymb,amsthm,units,stmaryrd}
\usepackage{qtree,bussproofs}

\usepackage{yfonts}
\usepackage{units,stmaryrd}
\usepackage[colorlinks]{hyperref}
\usepackage[usenames,dvipsnames]{color}
\hypersetup{
  colorlinks,
  citecolor=Violet,
  linkcolor=Red,
  urlcolor=Blue}
\usepackage{graphicx}
\usepackage[all]{xy}
\newtheorem{theorem}{Theorem}[section]

\newtheorem{definition}[theorem]{Definition}
\newtheorem{lemma}[theorem]{Lemma}
\newtheorem{remark}[theorem]{Remark}

\newtheorem{proposition}[theorem]{Proposition}

\usepackage{xspace}

\newcommand{\la}{\langle \,}
\newcommand{\ra}{\, \rangle}

\def\le{{\ell}}

\def\fmodels{\xymatrix{
\ar@{|=}[r]^{<\omega}&
}
}
\def\nmodels{\xymatrix{
\ar@{|=}[r]^{N}&
}
}
\def\<{\left <}

\def\>{\right >}

\DeclareSymbolFont{AMSb}{U}{msb}{m}{n}
\DeclareMathSymbol{\N}{\mathbin}{AMSb}{"4E}
\DeclareMathSymbol{\Z}{\mathbin}{AMSb}{"5A}
\DeclareMathSymbol{\R}{\mathbin}{AMSb}{"52}
\DeclareMathSymbol{\Q}{\mathbin}{AMSb}{"51}
\DeclareMathSymbol{\I}{\mathbin}{AMSb}{"49}

\newcommand{\bt}{\begin{theorem}}
\newcommand{\et}{\end{theorem}}
\newcommand{\bl}{\begin{lemma}}
\newcommand{\el}{\end{lemma}}

\newcommand\remove[1]{}

